\documentclass[11pt]{article}
\title{Characterization of the Fibonacci Cobweb Poset as oDAG}
\author{Ewa Krot \\
\\Institute of Computer Science, Bia{\l}ystok University\\
PL-15-887 Bia{\l}ystok, ul.Sosnowa 64, POLAND\\
e-mail: ewakrot@wp.pl, ewakrot@ii.uwb.edu.pl}

\usepackage{amsmath,amsthm}
\usepackage{graphicx}

\chardef\bslash=`\\ 
\newtheorem{thm}{Theorem}

\newtheorem{lemma}{Lemma}

\begin{document}
\maketitle
\begin{abstract}
The characterization of Fibonacci Cobweb poset $P$ as DAG and oDAG
is given. The {\em dim} 2 poset such that its Hasse diagram
coincide with digraf of $P$ is constructed.

\end{abstract}
\section{Fibonacci cobweb poset}
The Fibonacci cobweb poset $P$ has been invented by
A.K.Kwa\'sniewski in \cite{4,3,10} for the purpose of finding
combinatorial interpretation of fibonomial coefficients and
eventually their reccurence relation. 

In \cite{4} A. K. Kwa\'sniewski defined cobweb poset $P$ as
infinite labeled digraph oriented upwards as follows: Let us label
vertices of $P$ by pairs of coordinates: $\langle i,j \rangle \in
{\bf N_{0}}\times {\bf N_{0}}$, where the second coordinate is the
number of level in which the element of $P$ lies (here it is the
$j$-th level) and the first one is the number of this element in
his level (from left to the right), here $i$. Following \cite{4}
we shall refer to $\Phi_{s}$ as to the set of vertices (elements)
of the $s$-th level, i.e.:
$$\Phi_{s}=\left\{\langle j,s \rangle ,\;\;1\leq j \leq F_{s}\right\},\;\;\;s\in{\bf N}\cup\{0\},$$
where $\{F_{n}\}_{n\geq 0}$ stands for Fibonacci sequence.

Then $P$ is a labeled graph $P=\left(V,E\right)$ where
$$V=\bigcup_{p\geq 0}\Phi_{p},\;\;\;E=\left\{\langle \,\langle j,p\rangle,\langle q,p+1
\rangle\,\rangle\right\},\;\;1\leq j\leq F_{p},\;\;1\leq q\leq
F_{p+1}.$$

 We can now define the partial order relation on $P$ as follows:
let\\ $x=\langle s,t\rangle, y=\langle u,v\rangle $ be elements of
cobweb poset $P$. Then
$$ ( x \leq_{P} y) \Longleftrightarrow
 [(t<v)\vee (t=v \wedge s=u)].$$
 \vspace{2mm}
 \begin{center}
\includegraphics[width=80mm]{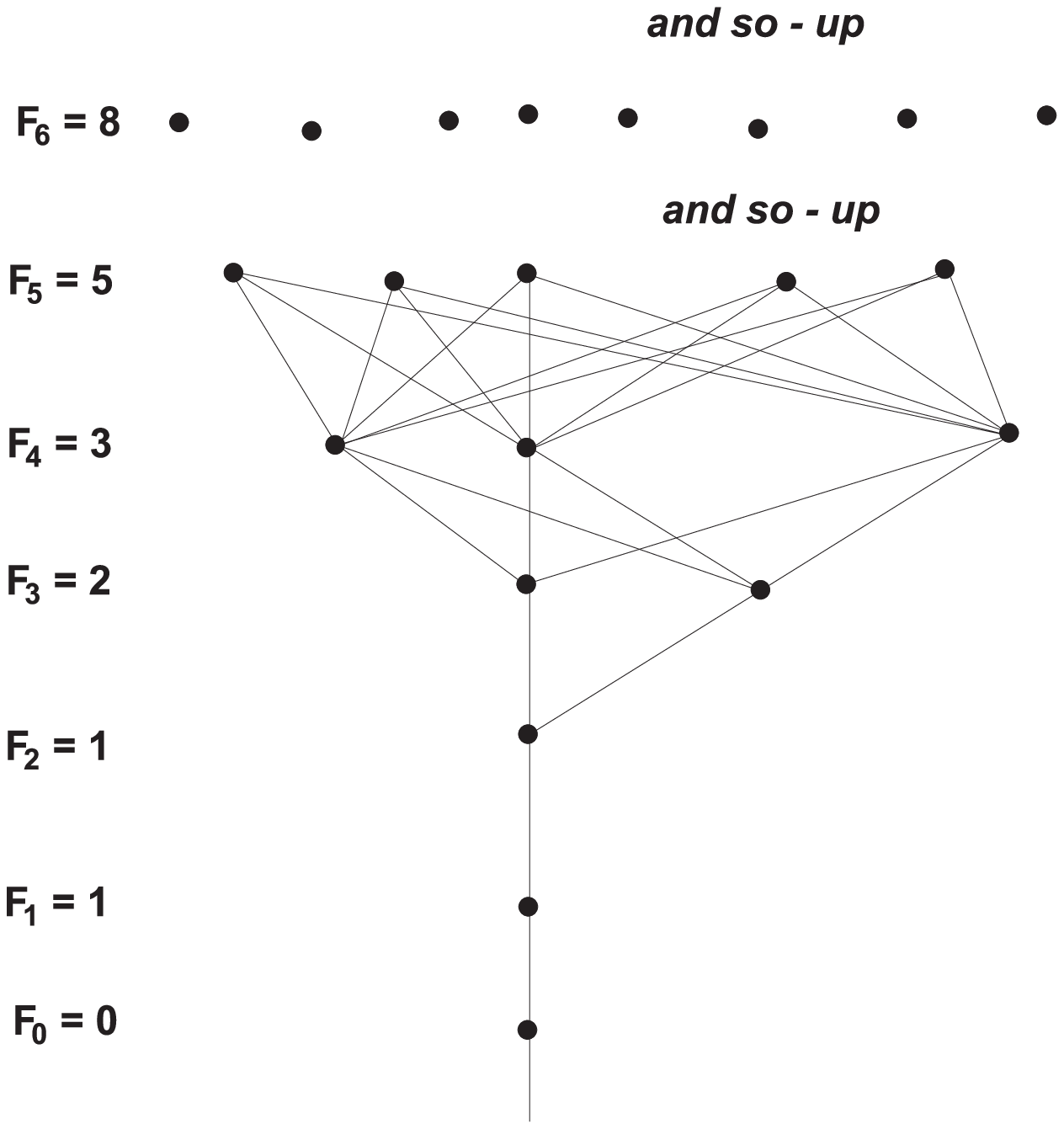}

\vspace{2mm}

\noindent {\small Fig.~1. The picture of the  Fibonacci "cobweb"
poset}
\end{center}
\section{DAG $\longrightarrow$ oDAG problem }
In \cite{p} A. D. Plotnikov considered the so called "DAG
$\longrightarrow$ oDAG problem". He determined condition when a
digraph $G$ may may be presented by the corresponding {\em dim } 2
poset $R$ and he established the algorithm for finding it.

Before citing Plotnikov's results lat us recall  (following
\cite{p})  some indispensable definitions.

If $P$ and $Q$ are partial orders on the same set $A$, $Q$ is said
to be an {\bf extension} of $P$ if $a\leq_{P} b$ implies
$a\leq_{Q} b$, for all $a,b\in A$. A poset $L$ is a {\bf chain},
or a {\bf linear order} if we have either $a\leq_{L} b$ or
$b\leq_{L} a$ for any $a,b\in A$. If $Q$ is a linear order then it
is a {\bf linear extension} of $P$.

The {\bf dimension} $dim\ R$ of $R$ being a partial order is the
least positive integer $s$ for which there exists a family $F=(L_1
,L_2 ,\ldots,L_s)$ of linear extensions of $R$ such that $R=
\bigcap_{i=1}^{s} L_{i}$. A family $F=(L_1,L_2,\ldots,L_s)$ of
linear orders on $A$ is called a {\bf realizer} of $R$ on $A$ if
\[
R=\bigcap_{i=1}^{s} L_{i}.
\]

We denote by $D_{n}$ the set of all acyclic directed $n$-vertex
graphs without loops and multiple edges. Each digraph ${\vec
G}=(V,{\vec E})\in D_{n}$ will be called {\bf DAG}.

A digraph ${\vec G}\in D_{n}$ will be called {\bf orderable
(oDAG)} if there exists are $dim\ 2$ poset such that its Hasse
diagram coincide with the digraph ${\vec G}$.

Let  ${\vec G}\in D_{n}$ be a digraph, which does not contain the
arc $(v_{i},v_{j})$ if there exists the directed path
$p(v_{i},v_{j})$ from the vertex $v_{i}$ into the vertex $v_{j}$
for any $v_{i}$, $v_{j}\in V$. Such digraph is called {\bf
regular}. Let $D\subset D_{n}$ is the set of all regular graphs.

Let there is a some regular digraph ${\vec G}=(V,E)\in D$, and let
the chain ${\vec X}$ has three elements $x_{i_{1}}$, $x_{i_{2}}$,
$x_{i_{3}}\in X$ such that $i_{1}<i_{2}<i_{3}$, and, in the
digraph ${\vec G}$, there are not paths $p(v_{i_{1}},v_{i_{2}})$,
$p(v_{i_{2}},v_{i_{3}})$ and there exists a path
$p(v_{i_{1}},v_{i_{3}})$. Such representation of graph vertices by
elements of the chain ${\vec X}$ is called the representation in
{\bf inadmissible form}. Otherwise, the chain ${\vec X}$ presets
the graph vertices in {\bf admissible form}.

 Plotnikov showed that:

\begin{lemma}{\em \cite{p}}\label{l1}
\label{l1} A digraph ${\vec G}\in D_{n}$ may be represented by a
$dim\ 2$ poset if:
\renewcommand{\labelenumi}{(\arabic{enumi})}
\begin{enumerate}
\item there exist two chains ${\vec X}$ and ${\vec Y}$, each of
which is a linear extension of ${\vec G}_{t}$; \item the chain
${\vec Y}$ is a modification of ${\vec X}$ with inversions, which
remove the ordered pairs of ${\vec X}$ that there do not exist in
${\vec G}$.
\end{enumerate}
\end{lemma}
Above lemma results in the algorithm for finding {\em dim} 2
representation of a given DAG (i.e. corresponding oDAG) while the
following theorem establishes the conditions for constructing it.
\begin{thm}{\em \cite{p}}\label{t1}
\label{t1} A digraph ${\vec G}=(V,{\vec E})\in D_{n}$ can be
represented by $dim\ 2$ poset iff it is regular and its vertices
can be presented by the chain ${\vec X}$ in admissible form.
\end{thm}
\section{Fibonacci cobweb poset as DAG and oDAG}
In this section we show that Fibonacci cobweb poset is a DAG and
it is orderable (oDAG).

Obviously, cobweb poset $P=(V, E)$ defined above is a DAG (it is
directed acyclic graph without loops and multiple edges). One can
also verify that it is regular. For two elements $\langle i,
n\rangle , \langle j,m\rangle \in V$ a directed path $p(\langle i,
n\rangle , \langle j,m\rangle)\notin E$ will esist iff $n<m+1$ but
then  $(\langle i, n\rangle , \langle j,m\rangle)\notin E$ i.e.
$P$ does not contain the edge $(\langle i, n\rangle , \langle
j,m\rangle)$.

It is also possible to verify that vertices of cobweb poset $P$
can be presented in admissible form by the chain ${\vec X}$ being
a linear extension of cobweb $P$ as follows:
\begin{multline*}{\vec X}=\Big(\langle 1,0\rangle,\langle 1,1\rangle ,
\langle 1,2\rangle, \langle 1,3\rangle, \langle 2,3\rangle,
\langle 1,4\rangle, \langle 2,4\rangle, \langle 3,4\rangle,\langle
1,5\rangle, \langle 2,5\rangle, \langle 3,5\rangle,\\\langle
4,5\rangle,\langle 5,5\rangle,...\Big),\end{multline*}
 where

$$ ( \langle s,t\rangle \leq_{{\vec X}} \langle u,v\rangle) \Longleftrightarrow
 [(s\leq u)\wedge (t\leq v)]$$
 for $1\leq s \leq F_{t},\; 1\leq u \leq F_{v},\;\;\;t, v \in {\bf N}\cup\{0\}.$

 Fibonacci cobweb poset $P$ satisfies the conditions of Theorem
 \ref{t1} so it is oDAG. To find the chain ${\vec Y}$ being
a linear extension of cobweb $P$ one uses Lemma \ref{l1} and
arrives at: \begin{multline*} {\vec Y}=\Big(\langle
1,0\rangle,\langle 1,1\rangle , \langle 1,2\rangle, \langle
2,3\rangle, \langle 1,3\rangle, \langle 3,4\rangle, \langle
2,4\rangle, \langle 1,4\rangle,\langle 5,5\rangle, \langle
4,5\rangle, \langle 3,5\rangle,\\\langle 2,5\rangle,\langle
1,5\rangle,...\Big),\end{multline*}
 where

$$ ( \langle s,t\rangle \leq_{{\vec Y}} \langle u,v\rangle) \Longleftrightarrow
 [(t < v)\vee (t=v \wedge s\geq u)]$$
 for $1\leq s \leq F_{t},\; 1\leq u \leq F_{v},\;\;\;t, v \in {\bf
 N}\cup\{0\}$   and finally
 $$ (P,\leq_{P})={\vec X}\cap{\vec Y}.$$
\section{Closing remark}
For any sequence $\{a_{n}\}$ of natural numbers one can define
corresponding cobweb poset as follows:
$$\Phi_{s}=\left\{\langle j,s \rangle ,\;\;1\leq j \leq a_{s}\right\},\;\;\;s\in{\bf N}\cup\{0\},$$
and  $cobP=\left(V,E\right)$ where
$$V=\bigcup_{p\geq 0}\Phi_{p},\;\;\;E=\left\{\langle \,\langle j,p\rangle,\langle q,p+1
\rangle\,\rangle\right\},\;\;1\leq j\leq a_{p},\;\;1\leq q\leq
a_{p+1}$$
 with the partial order relation on $cobP$ :
$$ ( x \leq_{P} y) \Longleftrightarrow
 [(t<v)\vee (t=v \wedge s=u)]$$
 for  $x=\langle s,t\rangle, y=\langle u,v\rangle $ being elements of
cobweb poset $cobP$. Similary as above one can show that the
family of cobweb posets consist of DAGs representable  by
corresponding {\em dim } 2 posets (i.e. of oDAGs).

 \end{document}